\documentclass[psamsfonts, reqno]{amsart}

\usepackage{amscd,amsmath,amssymb,amsfonts,verbatim}
\usepackage{times}
\usepackage[cmtip, all]{xy}
\usepackage{graphicx}
\usepackage[active]{srcltx}
\usepackage{color}
\usepackage{textcomp}

\definecolor{refkey}{gray}{.5}   
\definecolor{labelkey}{gray}{.5} 
\definecolor{Red}{rgb}{1,0,0}

\newcommand{\pf}{{\bf Proof : }}

\newtheorem{theo}{Theorem}[section]
\newtheorem{prop}[theo]{Proposition}
\newtheorem{lem}[theo]{Lemma}
\newtheorem{cor}[theo]{Corollary}

\theoremstyle{definition}

\newtheorem{notn}[theo]{Notation}
\newtheorem{rem}[theo]{Remark}

\newtheorem{defi}[theo]{Definition}

\newcommand{\Um}{\mbox{\rm Um}}		\newcommand{\SL}{\mbox{\rm SL}}
		
	\newcommand{\E}{\mbox{\rm E}}
\newcommand{\ESp}{\mbox{\rm ESp}}     \newcommand{\Sp}{\mbox{\rm Sp}}

\title{Unimodular rows over affine algebras over algebraic closure of a finite field}
\author{Sampat Sharma}
\newcommand{\Addresses}{{
  \bigskip
  \footnotesize

\medskip
  
 \textsc{Sampat Sharma, School of Basic Sciences, IIT Mandi, Mandi 175005 (H.P), India}\par\nopagebreak
  \textit{E-mail:} Sampat Sharma \texttt{<sampat@iitmandi.ac.in; sampat.iiserm@gmail.com>}

  \medskip

  }}

\begin{document}
\maketitle
\subjclass 2020 Mathematics Subject Classification:{13C10}

 ~ {Unimodular row, nice group structure}

 \begin{abstract}
In this article, we prove that if $R$ is an affine algebra of dimension $d\geq 4$ over $\overline{\mathbb{F}}_{p}$ and $1/(d-1)! \in R,$ then any unimodular row over $R$ of length $d$ can be mapped to a factorial row by elementary transformations.
 \end{abstract}
 
 \vskip0.50in

\section{Introduction} {\it All the rings are assumed to be commutative noetherian with unity $1$.}

One can say that the study of unimodular rows came in the picture with the Serre's question about projective modules. In \cite{jps}, Serre asked whether finitely generated projective modules over $k[X_1, \ldots, X_n]$ free, where $k$ is a field. This question of Serre was answered by Quillen \cite{quill} and Suslin \cite{4} independentaly. This is now known as Quillen-Suslin theorem. Due to the work of Grothendieck, one knows that every finitely generated projective module over $k[X_1, \ldots, X_n]$ is stably free. To a stably free module one can associate a unimodular row and vice versa. And the connection between these two objects is that a stably free module is free if and only if the corresponding unimodular row is completable to an invertible matrix. Therefore question about freeness of a projective module boils down to the question about completion of certain unimodular row. This shows that study of unimodular rows is quite important in classical $K$-theory.

 In \cite{swantowber}, R.G. Swan and J. Towber showed that if $(a^{2}, b, c)\in Um_{3}(R)$ then it can be completed to an invertible matrix over $R.$ This result of Swan and Towber was generalised by Suslin in \cite{sus2} who showed that if 
$(a_{0}^{r!}, a_{1}, \ldots, a_{r})\in Um_{r+1}(R)$ then it can be completed to an invertible matrix. Rao also studied the unimodular rows over $R[X]$ where $R$ is a local ring in \cite{invent}, \cite{trans}. In \cite{sampat}, author studied the unimodular rows over finitely generated rings over $\mathbb{Z}$ and proved that if $R$ is a finitely generated ring over $\mathbb{Z}$ of dimension $d, d\geq2, \frac{1}{d!}\in R$, then any unimodular row over $R[X]$ of length $d+1$ can be mapped to a factorial row by elementary transformations.

In \cite[Theorem 7.5]{frs}, Fasel, Rao and Swan proved that if $R$ is a normal affine algebra of dimension $d\geq 4$ over $\overline{\mathbb{F}}_{p}$ and $1/(d-1)! \in R,$ then any unimodular row over $R$ of length $d$ can be mapped to a factorial row by elementary transformations. In this article we remove their normality condition. In particular, we prove the following : 
\begin{theo} Let $R$ be an affine algebra of dimension $d\geq 4$ over $\overline{\mathbb{F}}_{p}$ and $1/(d-1)! \in R.$ Let 
$ v = (v_{1}, \ldots, v_{d})\in \Um_{d}(R).$ Then there exists $\varepsilon \in \E_{d}(R)$ such that 
$$(v_{1}, \ldots, v_{d})\varepsilon = (w_{1}, w_{2}, \ldots, w_{d}^{(d-1)!})$$ 
for some $(w_{1}, \ldots, w_{d})\in \Um_{d}(R).$ In particular, $\Um_{d}(R) = e_{1}\SL_{d}(R).$
\end{theo}

In this article, we consider the action of symplectic matrices on unimodular rows and prove the following : 
\begin{theo}  Let $R$ be an affine algebra of dimension $d$ over $\overline{\mathbb{F}}_{p}$ and $d \equiv 2~(\mbox{mod}~4).$ Let $v\in \Um_{d}(R).$ Assume that $\SL_{d+1}(R)\cap \E(R) = \E_{d+1}(R).$ Then $v$ is completable to a symplectic matrix.
\end{theo}

We also studied the problem of injective stability for affine algebras over $\overline{\mathbb{F}}_{p}$ and proved the following : 
\begin{theo}
\label{dropinjective}  Let $R$ be an affine algebra of dimension $d\geq 4$ over $\overline{\mathbb{F}}_{p}$. Let $\alpha \in \SL_{d}(R) \cap \E_{d+1}(R).$ Then $\alpha$ is isotopic to the identity matrix.
\end{theo}

\section{Preliminaries} 
\setcounter{theo}{0} In this section, we give some basic definitions and collect some known results which will be used later in the article.
\subsection{Basic definitions and notations}
\begin{defi}
 A row $v = (a_{0}, \ldots, a_{r})\in R^{r+1}$ is said to be unimodular if there is a $w 
 = (b_{0}, \ldots, b_{r})\in R^{r+1}$ with
$\langle v ,w\rangle = \Sigma_{i = 0}^{r} a_{i}b_{i} = 1$ and $Um_{r+1}(R)$ will denote
the set of unimodular rows (over $R$) of length 
$r+1$.
\end{defi}

\par 
 The elementary linear group 
$E_{r+1}(R)$ acts on the rows of length $r+1$ by right multiplication. Moreover, this action takes unimodular rows 
to unimodular
rows : $\frac{Um_{r+1}(R)}{E_{r+1}(R)}$ will denote the set of orbits of this action; and we shall denote by $[v]$ the 
equivalence class of a row $v$ under this equivalence relation. 

\begin{defi}
 A unimodular row $v\in Um_{n}(R)$ is said to be completable if there exists a matrix $\alpha \in GL_{n}(R)$ such that 
 $v = e_{1}\alpha.$
\end{defi}

\begin{notn}

Let $\psi_{1} = \begin{bmatrix}
                 0 & 1\\
                 -1 & 0\\
                \end{bmatrix},~~ \psi_{n} = \psi_{n-1} \perp \psi_{1}$ for $n > 1.$
\end{notn}

\begin{notn}
 Let $\sigma$ be the permutation of the natural numbers given by $\sigma(2i) = 2i-1$ and $\sigma(2i-1) = 2i$.
\end{notn}

\begin{defi}{\bf{Symplectic group}}~{$Sp_{2m}(R) :$} the group of all $2m \times 2m $ matrices 
$\{\alpha \in GL_{2m}(R)~\mid \alpha^{t}\psi_{m}\alpha = \psi_{m}\}$.
 
\end{defi}

\begin{defi}{\bf{Elementary symplectic group}}~{$ESp_{2m}(R)$}: We define for $ 1\leq i \neq j\leq 2m,~z\in R,$
$$
se_{ij}(z)=
\begin{cases}
I_{2m} + zE_{ij},~~~\textit{if}~ i = \sigma(j);\\
 I_{2m} + zE_{ij} - (-1)^{i+j}zE_{\sigma(j)\sigma(i)}, ~~~ \textit{if}~ i\neq \sigma(j).

\end{cases}
$$

It is easy to verify that all these matrices belong to $Sp_{2m}(R)$. We call them the elementary symplectic matrices
over $R$. The subgroup generated by them is called the elementary symplectic group and is denoted by $ESp_{2m}(R)$.
 
\end{defi}

\begin{defi}
     Let $I$ be an ideal of a ring $R.$ A unimodular row $v\in Um_{n}(R)$ which is congruent to $e_{1}$ modulo $I$ is 
     called unimodular relative to ideal $I.$ Set of unimodular rows relative to ideal $I$ will be denoted by $Um_{n}(R,I).$
    \end{defi}

Let $I$ be an ideal of a ring $R$. We shall denote by $GL_{n}(R,I)$ the kernel of the canonical mapping 
$GL_{n}(R)\longrightarrow GL_{n}(\frac{R}{I}).$ Let $SL_{n}(R,I)$ denotes the subgroup of
$GL_{n}(R,I)$ of elements of determinant $1$.

Let $I$ be an ideal of a ring $R$. We shall denote by $Sp_{2m}(R,I)$ the kernel of the canonical mapping 
$Sp_{2m}(R)\longrightarrow Sp_{2m}(\frac{R}{I}).$

\begin{defi}${\bf{The~ relative~ groups~ E_{n}(I),~ E_{n}(R,I):}}$
 Let $I$ be an ideal of $R$. The elementary group $E_{n}(I)$ is the subgroup of $E_{n}(R)$ generated
 as a group by the elements
 $e_{ij}(x),~x\in I,~1\leq i\neq j\leq n.$\\
  The relative elementary group $E_{n}(R,I)$ is the normal closure of $E_{n}(I)$ in $E_{n}(R)$.
  \end{defi}
  
  \begin{defi}${\bf{The~ relative~ groups~ ESp_{2m}(I),~ ESp_{2m}(R,I):}}$
 Let $I$ be an ideal of $R$. The elementary symplectic group $ESp_{2m}(I)$ is the subgroup of $ESp_{2m}(R)$ 
 generated as a group by the 
 elements $se_{ij}(x),~x\in I,~1\leq i\neq j\leq 2m.$
 \par
  The relative elementary symplectic group $ESp_{2m}(R,I)$ is the normal closure of $ESp_{2m}(I)$ in $ESp_{2m}(R)$.
  \end{defi}
  
\subsection{The Suslin matrices}
    \setcounter {theo}{0}
    First recall the Suslin matrix $S_{r}(v,w)$. These were defined by Suslin in 
$(${\cite[Section 5] {sus2}}$).$ We recall his inductive process : Let $v = (a_{0}, v_{1}), w = (b_{0}, w_{1})$, where 
$v_{1}, w_{1} \in M_{1,r}(R)$. Set $S_{0}(v,w) = a_{0}$ and set 
$$S_{r}(v,w) = \begin{bmatrix}
                 a_{0}I_{2^{r-1}} & S_{r-1}(v_{1},w_{1})\\
                 -S_{r-1}(w_{1}, v_{1})^{T} & b_{0}I_{2^{r-1}}\\
                \end{bmatrix}.$$

\par 
The process is reversible and given a Suslin matrix $S_{r}(v,w)$ one can recover the associated rows $v,w$, i.e. the pair 
$(v,w)$. 

\par 
 In \cite{suskcohomo}, Suslin proves that if $\langle v ,w\rangle = v\cdot w^{T} = 1$, then by row and column operations
one can reduce 
$S_{r}(v,w)$ to a matrix $\beta_{r}(v,w)$ of size $r+1$ whose first row is $(a_{0},a_{1},a_{2}^{2}\cdots,a_{r}^{r})$. 
This in particular proves that rows of such type can be completed to an invertible matrix of determinant $1.$ We call 
$\beta_{r}(v,w)$ to be a compressed Suslin matrix. 

\subsection{Some assorted results}
 \setcounter {theo}{0}
We state two results from \cite[Corollary 17.3, Corollary 18.1, Theorem 18.2]{7}.

\begin{theo}
\label{transit} Let $R$ be an affine $C$-algebra of dimension $d\geq 2$, where $C$ is either a subfield $F$ of  $\overline{\mathbb{F}}_{p}$ or $C=\mathbb Z$. Then 
\begin{itemize}
\item If $d =2,$ then ${\E}_{3}(R)$ acts transitively on $\Um_{3}(R).$
\item If $d\geq 3,$ then $\mbox{sr}(R) \leq d.$ As a consequence, ${\E}_{d+1}(R)$ acts transitively on $\Um_{d+1}(R).$
\end{itemize}
\end{theo}

We will state the part of Swan's Bertini needed in our proof \cite[Theorems 1.3, 1.4]{swanbertini} (see also \cite[Theorem 2.3]{murthy}).

\begin{theo}\label{sb}
Let $R$ be a geometrically reduced affine ring of dimension $d$ over an infinite field $k$ and $(a,a_1,\ldots,a_r)\in \Um_{r+1}(R)$. Then there exist $b_1,\ldots,b_r\in R$ such that if $a'=a+a_1b_1+\ldots+a_rb_r$, then $R/(a')$ is a reduced affine algebra  of dimension $d-1$ which is smooth at all smooth points of $R$. 
\end{theo}

We note results of Fasel, Rao and Swan in \cite[]{frs}. 
\begin{theo}\label{frslabel} Let $R$ be a smooth affine algebra of dimension $d\geq 3$ over an algebraically closed field $k$ and $\mbox{gcd}(d-1)!, \mbox{char}(k)) = 1.$ Let $v = (v_{1}, \ldots, v_{d})\in \Um_{d}(R).$ Then there exists $\varepsilon \in E_{d}(R)$ such that 
$$v\varepsilon = (w_{1}, w_{2}, \ldots, w_{d}^{(d-1)!})$$ 
for some $ (w_{1}, w_{2}, \ldots, w_{d})\in \Um_{d}(R).$
\end{theo} 

\begin{cor}\cite[Corollary 7.7]{frs} 
Let $R$ be a smooth affine algebra of dimension $d\geq 3$ over an algebraically closed field $k.$ Assume that $d! k = k.$ Then $SL_{d}(R) \cap E_{d+1}(R) = E_{d}(R).$
\end{cor}

\begin{cor}\cite[Corollary 7.8]{frs}
Let $R$ be a smooth affine algebra over the algebraic closure $k$ of a finite field of dimension $d\geq 3.$ Assume that $d! k = k.$ Then, the natural map $\frac{SL_{d}(R)}{E_{d}(R)} \to SK_{1}(R)$ is an isomorphism.
\end{cor}

\begin{cor}\cite[Corollary 7.9]{frs} Let $R$ be a smooth affine algebra of dimension $d\leq 4$ over an algebraically closed field $k.$ Assume that $d! k = k.$ Then, the Vaserstein symbol $V: Um_{3}(R)/E_{3}(R) \to W_{E}(R)$ is an isomorphism.
\end{cor}

\par Next we note a result of Vaserstein  {\cite[Lemma 5.5]{7}}.
\begin{lem}
\label{2.11}
  Let $R$ be a ring. Then for any $m\geq 1,$
 ${\E}_{2m}(R)e_{1} = ({\Sp}_{2m}(R) \cap {\E}_{2m}(R))e_{1}.$
\end{lem}
 
\begin{rem}
\label{2.12}
 It was observed in {\cite[Lemma 2.13]{pr}} that Vaserstein's proof actually shows that  
 ${\E}_{2m}(R)e_{1} = {\ESp}_{2m}(R)e_{1}.$ 
\end{rem}

We recall a result of Basu, Chattopadhyay and Rao {\cite[Lemma 2.13]{bcr}}.
\begin{lem}
\label{bcrlem} Let $R$ be a ring essentially of dimension $d\geq 1$ and $S_{r}(v,w)\in \SL_{2^{r}(R, I)}, 2^{r}\geq d+2$ for $v, w\in \Um_{r+1}(R, I)$ satisfying $\langle v, w\rangle = 1.$ Let us assume that $r\equiv 1~(\mbox{mod}~4).$ Then there exists $\varepsilon_{J_{r}} \in \ESp_{J_{r}}(R, I)$ such that $S_{r}(v,w)\varepsilon_{J_{r}} = (I_{2^{r}-k}\perp \gamma)\varepsilon$ for some $\gamma \in \Sp_{k}(R,I)$ and $\varepsilon \in \E_{2^{r}}(R,  I).$ Here $k = d+1$ if $d$ is odd and $k = d$ if 
$d$ is even.
\end{lem}

We note the following two results of Vaserstein in \cite{vas}:

\begin{cor}
If the ring $\frac{R}{rad{R}}$ decomposes into product of any number of matrix rings over nonassociative division rings (For example, $R$ is a finite dimension algebra over a finite field), then stable range of $R = 1.$ 
\end{cor}

\begin{cor}
The stable range of $R^{1} = \mbox{max}(2, st. r. (R))$ where $R^{1}$ is the ring obtained by formal adjunction of unity to $R.$ Note that $\frac{R^{1}}{R} = \mathbb{Z}.$
\end{cor}
\section{Completion to a special linear matrix}
\setcounter{theo}{0}
In this section we prove that unimodular rows of length $d$ over $R$ can be elementarily mapped to a factorial row when R is an affine algebra of dimension $d\geq 4$ over $\overline{\mathbb{F}}_{p}$ and $1/(d-1)! \in R.$ We also prove this result in the relative case.
\begin{theo}\label{completion} Let $R$ be an affine algebra of dimension $d\geq 4$ over $\overline{\mathbb{F}}_{p}$ and $1/(d-1)! \in R.$ Let 
$ v = (v_{1}, \ldots, v_{d})\in \Um_{d}(R).$ Then there exists $\varepsilon \in \E_{d}(R)$ such that 
$$(v_{1}, \ldots, v_{d})\varepsilon = (w_{1}, w_{2}, \ldots, w_{d}^{(d-1)!})$$ 
for some $(w_{1}, \ldots, w_{d})\in \Um_{d}(R).$ In particular, $\Um_{d}(R) = e_{1}\SL_{d}(R).$
\end{theo}
${\pf}$ In view of {\cite[Remark 1.4.3]{invent}},we may assume that $R$ is reduced. Let $J$ be the ideal defining singular locus of $R$. Since $R$ is a reduced ring, $\mbox{ht}(J)\geq 1.$ Let $\overline{R} = R/J.$ Thus $\mbox{dim}(\overline{R})\leq d-1$ and $ \overline{v} = (\overline{v}_{1}, \ldots, \overline{v}_{d})\in \Um_{d}(\overline{R}).$ In view of Theorem \ref{transit}, $\mbox{sr}(\overline{R}) \leq d-1.$ Thus there exist $\overline{\lambda}_{2}, \ldots, \overline{\lambda}_{d}\in \overline{R}$ such that 
$(\overline{v}_{2}', \ldots, \overline{v}_{d}')\in \Um_{d-1}(\overline{R}),$
where 
$\overline{v}_{i}' = \overline{v}_{i} + \overline{\lambda_{i}v_{i}}, 2\leq i\leq d.$ Note that $[(\overline{v}_{1}, \overline{v}_{2}, \ldots, \overline{v}_{d})] = [(\overline{v}_{1}, \overline{v}_{2}', \ldots, \overline{v}_{d}')].$

Since $(\overline{v}_{2}', \ldots, \overline{v}_{d}')\in \Um_{d-1}(\overline{R}),$ we can elementarily transform $(\overline{v}_{1}, \overline{v}_{2}', \ldots, \overline{v}_{d}')$ to $(1, 0, \ldots, 0). $ Thus we may assume that $(v_{1}, \ldots, v_{d}) \equiv e_{1}~\mbox{mod}~J.$ Now in view of the Theorem \ref{sb}, we add multiples of $v_{2}, \ldots, v_{d}$ to $v_{1}$ to transform it to $v_{1}'$ such that $R/(v_{1}')$ is smooth at the smooth points of $R.$ Since $v_{1}' \equiv 1~\mbox{mod}~J$, 
$R/(v_{1}')$ is smooth of dimension $d-1.$

Thus in view of the Theorem \ref{frslabel}, there exists $\overline{\varepsilon}\in \E_{d-1}(\overline{R})$ such that 
$(\overline{v}_{2}, \ldots, \overline{v}_{d})\overline{\varepsilon} = (\overline{w}_{2}, \ldots, \overline{w}_{d}^{(d-1)!}).$ Let $\varepsilon$ be a lift of $\overline{\varepsilon}$ and upon making appropriate elementary transformations we get 
$$(v_{1}, v_{2},\ldots, v_{d})\varepsilon = (w_{1}, w_{2}, \ldots, w_{d}^{(d-1)!}).$$
$\hfill \square$

\begin{lem}
\label{relelementary} Let $R$ be an affine algebra of dimension $d\geq 4$ over $\overline{\mathbb{F}}_{p}$ and let $v\in \Um_{d}(R, I)$ for some ideal $I \subset R.$ Let $J$ denotes the ideal defining singular locus of $R$. Then there exists $\varepsilon \in \E_{d}(R, I)$ such that 
$v\varepsilon = (u_{1}, u_{2}, \ldots, u_{d})$ with $u_{1} \equiv 1~\mbox{mod}~(I \cap J).$
\end{lem}
${\pf}$ Let $v = (v_{1}, v_{2},\ldots, v_{d})\in \Um_{d}(R, I)$ and $\overline{R} = R/J.$ We have $\overline{v} = (\overline{v_{1}}, \overline{v_{2}},\ldots, \overline{v_{d}})\in \Um_{d}(\overline{R}, \overline{I}).$ In view of the Theorem \ref{transit},  $\mbox{sr}(\overline{R}) \leq d-1.$  Thus upon adding multiples of $\overline{v_{d}}$ to $\overline{v_{1}}, \overline{v_{2}},\ldots, \overline{v_{d-1}}$,  we may assume that $(\overline{v_{1}}, \overline{v_{2}},\ldots, \overline{v_{d-1}})\in \Um_{d-1}(\overline{R}, \overline{I}).$ Let $(\overline{w_{1}}, \overline{w_{2}},\ldots, \overline{w_{d-1}})\in \Um_{d-1}(\overline{R}, \overline{I})$ be such that 
$\sum_{i=1}^{d-1} \overline{v_{i}} \overline{w_{i}} = \overline{1}.$ Let $$\overline{\varepsilon_{1}} = \begin{bmatrix}
                 \overline{1} &  \overline{0} & \ldots & (1-\overline{v_{d}})\overline{w_{1}}\\
                  \overline{0} & \overline{1}& \ldots & (1-\overline{v_{d}})\overline{w_{2}}\\
\vdots & & & \vdots \\
\overline{0} & \overline{0}& \ldots & (1-\overline{v_{d}})\overline{w_{d-1}}\\
\overline{0} & \overline{0}& \ldots & \overline{1}\\
                \end{bmatrix} \in \E_{d}(\overline{R}).$$ Then $(\overline{v_{1}}, \overline{v_{2}},\ldots, \overline{v_{d}})\overline{\varepsilon_1} = (\overline{v_{1}}, \overline{v_{2}},\ldots, \overline{v_{d-1}}, \overline{1}).$ There exists $\overline{\varepsilon_2}
\in\E_{d}(\overline{R}, \overline{I})$ such that $(\overline{v_{1}}, \overline{v_{2}},\ldots, \overline{v_{d-1}}, \overline{1})\overline{\varepsilon_2} = (\overline{1}, \overline{0},\ldots, \overline{0}, \overline{1}).$ Now $(\overline{1}, \overline{0},\ldots, \overline{0}, \overline{1})\overline{\varepsilon_1}^{-1} = (\overline{1}, \overline{0},\ldots, \overline{0}, \overline{1}-(1-\overline{v_{d}})\overline{w_1}).$ Since 
$1 -(1-\overline{v_{d}})\overline{w_1}) = \overline{0}~\mbox{mod}~{\overline{I}},$ there exists $\overline{\varepsilon_3}\in \E_{d}(\overline{R}, \overline{I})$ such that   $ (\overline{1}, \overline{0},\ldots, \overline{0}, 1-(1-\overline{v_{d}})\overline{w_1})\overline{\varepsilon_3} = (\overline{1}, \overline{0},\ldots, \overline{0}, \overline{0}).$  Let $\overline{\varepsilon} = \overline{\varepsilon_1}\overline{\varepsilon_2}\overline{\varepsilon_1}^{-1} \overline{\varepsilon_3} \in \E_{d}(\overline{R}, \overline{I}).$ Let  $\varepsilon \in \E_{d}(R, I)$ be a lift of $\overline{\varepsilon}.$ Then we have $v\varepsilon = (u_{1}, u_{2}, \ldots, u_{d})$ with $u_{1} \equiv 1~\mbox{mod}~(I \cap J).$
$\hfill \square$


\begin{theo}
\label{completionrelative} Let $R$ be an affine algebra of dimension $d\geq 4$ over $\overline{\mathbb{F}}_{p}$ and let $v\in \Um_{d}(R, I)$ for some ideal $I \subset R.$ Then there exists $\varepsilon \in \E_{d}(R, I)$ such that 
$$v\varepsilon = (w_{1}, \ldots, w_{d}^{(d-1)!})$$ 
for some $(w_{1}, \ldots, w_{d})\in \Um_{d}(R, I).$
\end{theo}

${\pf}$ Let $v = (v_{1}, \ldots, v_{d})\in \Um_{d}(R,I)$ and $J$ be the ideal defining singular locus of $R.$ In view of the Lemma \ref{relelementary}, we may assume that $v_{1} \equiv 1~\mbox{mod}~(I\cap J).$ Let $v_{1} = 1 - \lambda$ for some $\lambda \in I\cap J.$ Note that $ (v_{1}, \ldots, v_{d}) \underset{\E_{d}(R, I)}{\sim}  (v_{1}, \lambda v_{2},\ldots, \lambda v_{d}).$ Now in view of the Theorem \ref{sb}, we add multiples of $\lambda v_{2}, \ldots, \lambda v_{d}$ to $v_{1},$ to transform it to $v_{1}'$ such that $R/(v_{1}')$ is smooth outside the singular set of $R$. Since $v_{1}' \equiv 1~\mbox{mod}~ J$, $R/(v_{1}')$ is a smooth affine algebra of dimension $d-1.$

Let $v_{1}' = 1 - \eta$ for some $\eta \in I\cap J.$ Note that in $\overline{R} = R/(v_{1}'),$ we have $\overline{\eta} = \overline{1}.$ In view of the Theorem \ref{frslabel}, there exists $\overline{\varepsilon}\in \E_{d-1}(\overline{R}) (= \E_{d-1}(\overline{\eta}))$ such that $(\overline{ \lambda v_{2}}, \ldots, \overline{\lambda v_{d}})\overline{\varepsilon} = 
(\overline{w_{2}}, \ldots, \overline{w_{d}}^{(d-1)!}).$ Let $\varepsilon \in \E_{d-1}(\eta) \subset \E_{d-1}(I)$ be a lift of $\overline{\varepsilon}.$ Therefore $({ \lambda v_{2}}, \ldots, {\lambda v_{d}}){\varepsilon} = 
({w_{2}}, \ldots, {w_{d}}^{(d-1)!})~\mbox{mod}~(Iv_{1}').$ Thus, we have 
\begin{align*} 
(v_{1}, \ldots, v_{d}) &\underset{\E_{d}(R, I)}{\sim} (v_{1}', \lambda v_{2},\ldots, \lambda v_{d})\\
& \underset{\E_{d}(R, I)}{\sim} (w_{1},{w_{2}}, \ldots, {w_{d}}^{(d-1)!}),~~~~~~\mbox{where}~w_{1} = v_{1}'.
\end{align*}
$\hfill \square$

Now we consider the problem of injective stability and prove that stably elementary matrices of size $d$ are isotopic to the identity. To prove this we need the following proposition.

\begin{prop}
\label{completionrelatives} Let $R$ be an affine algebra of dimension $d\geq 4$ over $\overline{\mathbb{F}}_{p}$. Let $v\in \Um_{d}(R)$ be such that $v \equiv e_{1}~\mbox{mod}~(s)$ for some $s\in R.$ Then there exists $\sigma \in \SL_{d}(R)$ with 
$\sigma \equiv \mbox{I}_{d}~\mbox{mod}~(s)$ such that $v = e_{1}\sigma.$
\end{prop}
${\pf}$ Let $S = R[X]/(X^{2}-sX).$ Then $S$ is an affine algebra of dimension $d\geq 4$ over $\overline{\mathbb{F}}_{p}$. Since $v \equiv e_{1}~\mbox{mod}~(s),$ let $v = e_{1} + sw$ for some $w\in R^{d}.$ Let $u(X) = e_{1}+Xw.$ \\
{\bf{Claim :}} $u(X) \in \Um_{d}(S).$\\
\noindent
{\bf{Proof of the claim :}} Since $v = e_{1} + sw,$ there exist $w_{1}', \ldots, w_{d}'\in R $ such that 
$$(1+sw_{1})w_{1}' + sw_{2}w_{2}' + \cdots + sw_{d}w_{d}' = 1.$$   Suppose to the contrary that $u(X) \notin \Um_{d}(S).$ Thus $\langle 1 + Xw_{1}, Xw_{2}, \ldots, Xw_{d}\rangle \subset \mathfrak{p}$ for some $\mathfrak{p} \in \mbox{Spec}(S).$
Therefore \begin{align*}
(1+Xw_{1})w_{1}' + \cdots + Xw_{d}w_{d}'  & = w_{1}' + X(w_{1}w_{1}' + \cdots w_{d}w_{d}')\\
& = 1 + (X-s)(w_{1}w_{1}' + \cdots w_{d}w_{d}') \in \mathfrak{p}.
\end{align*}
Since $X(X-s) = 0$ in $S$, $X\in \mathfrak{p}.$ As $1+Xw_{1} \in \mathfrak{p},$ $1\in \mathfrak{p}$ which is not possible. Thus $u(X) \in \Um_{d}(S).$
\par Now $u(s) = v$ and $u(0) = e_{1}.$ By Theorem \ref{completion}, there exists $\alpha(X) \in \SL_{d}(S)$ such that $u(X)  = e_{1}\alpha(X).$ Upon taking $\sigma = \alpha(0)^{-1}\alpha(s),$ we have $v = e_{1}\sigma$ and $\sigma \equiv \mbox{I}_{d}~\mbox{mod}~(s).$
$\hfill \square$

\begin{theo}
\label{dropinjective}  Let $R$ be an affine algebra of dimension $d\geq 4$ over $\overline{\mathbb{F}}_{p}$. Let $\alpha \in \SL_{d}(R) \cap \E_{d+1}(R).$ Then $\alpha$ is isotopic to the identity matrix.
\end{theo}
${\pf}$ Since $\alpha \in \SL_{d}(R)\cap \E_{d+1}(R)$, there exists a $\beta(T) \in \SL_{d+1}(R[T])$ such that 
$$\beta(0) = I_{d+1}~\mbox{and}~\beta(1) = 1\perp \alpha.$$
Let $A = R[T], s = T^{2} - T$ and $v = e_{1}\beta(T).$ By Proposition \ref{completionrelatives}, there exists $\gamma(T) \in \SL_{d+1}(R[T], (s))$ such that $v = e_{1}\gamma(T).$ Therefore $e_{1}\beta(T)\gamma(T)^{-1} = e_{1}.$ Thus 
$$\beta(T)\gamma(T)^{-1} = (1\perp \eta(T))\Pi_{i=2}^{d+2}E_{i1}(r_{i}),$$ 
where $\eta(T)\in \SL_{d+1}(R[T])$ and $r_{i}\in R[T].$ Now, 
$$\beta(0)\gamma(0)^{-1} = I_{d+1} = 1\perp \eta(0),$$
and 
$$\beta(1)\gamma(1)^{-1} = (1\perp \alpha)I_{d+1} = 1\perp \eta(1).$$
Thus $\eta(T)$ is an isotopy from $\alpha$ to $I_{d+1}.$
$\hfill \square$

\section{Completion to a symplectic matrix}
\setcounter{theo}{0}
In this section, we prove that symplectic matrices of size $d$ acts transitively on unimodular rows of length $d$ when $d \equiv 2~(\mbox{mod}~4).$
\begin{lem}  Let $R$ be an affine algebra of dimension $d$ over $\overline{\mathbb{F}}_{p}$ and $d \equiv 2~(\mbox{mod}~4).$ Let $v\in \Um_{d}(R).$ Assume that $\SL_{d+1}(R)\cap \E(R) = \E_{d+1}(R).$ Then $v$ is completable to a symplectic matrix.
\end{lem}
${\pf}$ Let $d = 2$ and $v = (v_{1}, v_{2})\in \Um_{2}(R).$ Let $w_{1}, w_{2}\in R$ be such that $v_{1}w_{1} + v_{2}w_{2} = 1.$ Consider 
$$\alpha =  \begin{bmatrix}
                 v_{1} &  v_{2} \\
                  -w_{2} & w_{1}\\
                \end{bmatrix} \in \SL_{2}(R) = \Sp_{2}(R).$$ 
Thus $v$ is completable to a symplectic matrix.

Thus we may assume that $d\geq 6.$ In view of Theorem \ref{completion}, there exists $\varepsilon \in \E_{d}(R)$ such that $v\varepsilon = \chi_{d!}(v_{1})$ for some $v_{1}\in \Um_{d}(R).$  By Remark \ref{2.12}, there exists $\varepsilon_{1}\in \ESp_{d}(R)$ such that $v\varepsilon_{1} = \chi_{d!}(v_{1}).$ Let $w_{1}\in \Um_{d}(R)$ be such that $\langle v_{1}, w_{1}\rangle = 1.$ By Lemma \ref{bcrlem}, we have $\varepsilon_{J_{d-1}} \in \ESp_{J_{d-1}}(R)$ such that 
$$S_{d-1}(v,w)\varepsilon_{J_{d-1}} = (I_{2^{d-1}-d}\perp \gamma)\varepsilon_{2}$$
for some $\gamma \in \Sp_{d}(R)$ and $\varepsilon_{2}\in \E_{2^{d-1}}(R).$ Therefore $S_{d-1}(v,w)$ and $\gamma$ are stably elementary equivalent. One knows that  $S_{d-1}(v_{1},w_{1})$ and $\beta_{d-1}(v_{1},w_{1})$ are stably elementary equivalent, we have 
$\beta_{d-1}(v_{1},w_1)\gamma^{-1}\in \SL_{d}(R)\cap \E_{d+1}(R).$ By {\cite[Corollary 3.9]{kesharisharma2}}, first row of an $1$-stably elementary matrix is elementarily completable. Thus $e_{1}\beta_{d-1}(v_1,w_1)\gamma^{-1} = e_{1}\delta'$ for some $\delta'\in \E_{d}(R).$ Now by Remark \ref{2.12}, $e_{1}\delta' = e_{1}\delta$ for some $\delta \in \ESp_{d}(R).$ Therefore $$e_{1}\beta_{d-1}(v_1,w_1)\gamma^{-1} = e_{1}\delta.$$ 
Thus $\chi_{d!}(v_1) = e_{1}\beta_{d-1}(v_1,w_1) = e_{1}\delta\gamma.$ Since $v\varepsilon_{1} = \chi_{d!}(v_1),$ we have $v = e_{1}\delta\gamma\varepsilon_{1}^{-1}$ and $\delta\gamma\varepsilon_{1}^{-1}\in \Sp_{d}(R).$
$\hfill \square$

\begin{prop}
\label{propsymplectic} Let $R$ be an affine algebra of dimension $d$ over $\overline{\mathbb{F}}_{p}$ and $d \equiv 2~\mbox{mod}~4.$ Let $v\in \Um_{d}(R)$ be such that $v\equiv e_{1}~\mbox{mod}~(s)$ for some $s\in R.$ Assume that $\SL_{d+1}(R)\cap \E(R) = \E_{d+1}(R).$ Then there exists $\sigma \in \Sp_{d}(R)$ be such that $\sigma \equiv \mbox{I}_{d}~(\mbox{mod}~s)$ and $v = e_1\sigma.$
\end{prop}
${\pf}$ The arguments of this proof are very similar to the proof of Proposition \ref{completionrelatives}.
$\hfill \square$

The next result shows that a stably elementary symplectic matrix over $R$ is isotopic to the identity matrix under certain conditions when $R$ is an affine algebra of dimension $d$ over $\overline{\mathbb{F}}_{p}.$ 

\begin{theo}
\label{dropsymplectic}  Let $R$ be an affine algebra of dimension $d$ over $\overline{\mathbb{F}}_{p}.$ Assume that $1/(d+1)! \in \overline{\mathbb{F}}_{p}$ and $d\equiv 1~(\mbox{mod}~4).$ Further assume that $\SL_{d+1}(R)\cap \E(R) = \E_{d+1}(R).$ Let $\sigma \in \Sp_{d-1}(R)\cap \ESp_{d+1}(R).$ Then $\sigma$ is isotopic to the identity.
\end{theo}
${\pf}$ Since $\sigma \in \Sp_{d-1}(R)\cap \ESp_{d+1}(R)$, there exists $\delta(T)\in \ESp_{d+1}(R[T])$ be such that $\delta(1) = I_{2}\perp \sigma, \delta(0) = \mbox{I}_{d+1}.$ Let $v(T) = e_{1}\delta(T)\in \Um_{d+1}(R[T], (T^{2}-T)).$ By Proposition \ref{propsymplectic}, $v(T) = e_{1}\alpha(T)$ for some $\alpha(T)\in \Sp_{d+1}(R[T])$ and $\alpha(T)\equiv \mbox{I}_{d+1}~\mbox{mod}~(T^{2}-T).$ Thus $e_{1}\delta(T) = e_{1}\alpha(T).$ Therefore one has 
$e_{1}\delta(T)\alpha(T)^{-1} = e_{1}.$ Now note that $\delta(T)\alpha(T)^{-1} $ is an isotopy of $\mbox{I}_{2}\perp \sigma$ to $\mbox{I}_{d+1}.$ The matrix $\delta(T)\alpha(T)^{-1} $ looks like 
$$ \delta(T)\alpha(T)^{-1} =  \begin{bmatrix}
                 1 &  0 & 0 \\
                  \ast & 1 & \ast\\
\ast & 0 & \eta(T)\\
                \end{bmatrix}$$
for some $\eta(T) \in \Sp_{d-1}(R[T], (T^{2}-T)).$ Now note that $\eta(0) = \mbox{I}_{d-1}$ and $\eta(1) = \sigma.$ Thus $\sigma$ is isotopic to the identity.
$\hfill \square$

\Addresses

\end{document}